\documentclass[12pt]{article}
\usepackage{amsmath, amssymb, amsthm}
\usepackage{mathrsfs}
\usepackage{hyperref}

\numberwithin{equation}{section}

\newtheorem{theorem}{Theorem}[section]
\newtheorem{lemma}[theorem]{Lemma}
\newtheorem{definition}[theorem]{Definition}
\newtheorem{proposition}[theorem]{Proposition}

\newtheorem{corollary}[theorem]{Corollary}

\newcommand{\N}{\mathbb{N}}

\newcommand{\diam}{\operatorname{diam}}

\title{Fixed Point Theorems for Relaxed Asymptotic Contractions via Two Quasi-Metrics}
\author{Jie Shi\\ Department of Mathematics and Statistics, Hubei Engineering University\\ \texttt{shijie@whu.edu.cn}}
\date{April 17, 2026}

\begin{document}

\maketitle

\begin{abstract}
We introduce a new class of asymptotic contractions that employs two quasi-metrics defined directly in terms of the underlying mapping. The contraction condition compares these two quantities via a sequence of bounding functions that converge locally uniformly to a Boyd–Wong function. This framework relaxes the hypotheses of Kirk's asymptotic fixed point theorem and strictly contains it as a special case. Assuming only the continuity of the map and the boundedness of some orbit in a complete metric space, we prove both the existence and uniqueness of a fixed point, along with the convergence of all iterates to that point.
\end{abstract}

\noindent\textbf{Keywords:} fixed point; asymptotic contraction; lower quasi-metric; upper quasi-metric; adaptive switching; local uniform convergence; Boyd–Wong condition

\section{Introduction}

The evolution of metric fixed point theory has been propelled by a persistent weakening of contractive hypotheses. In 1912, Brouwer \cite{Brouwer1912} established the topological fixed point principle for compact convex sets in Euclidean spaces. A decade later, Banach \cite{Banach1922} provided the fundamental contraction mapping theorem, which requires a uniform Lipschitz constant $\alpha<1$ and guarantees both existence and uniqueness of fixed points in complete metric spaces:
\begin{equation}
d(Tx,Ty) \leq \alpha\, d(x,y) \qquad (\alpha<1). \label{eq:banach}
\end{equation}
Schauder \cite{Schauder1930} subsequently extended fixed point theory to infinite-dimensional Banach spaces by means of compactness arguments.

The study of nonexpansive mappings, characterized by the condition
\begin{equation}
d(Tx,Ty)\leq d(x,y) \quad \text{for all } x,y, \label{eq:nonexpansive}
\end{equation}
gained prominence in the 1960s and revealed that additional geometric assumptions are indispensable for ensuring the existence of fixed points. Rakotch \cite{Rakotch1962} introduced nonlinear contractive conditions of the form
\begin{equation}
d(Tx,Ty)\leq \alpha(d(x,y))\,d(x,y) \label{eq:rakotch}
\end{equation}
with $\alpha$ decreasing and $\alpha(t)<1$ for $t>0$. Shortly thereafter, Edelstein \cite{Edelstein1962} proved the uniqueness of fixed points for strictly nonexpansive maps on compact metric spaces. The celebrated Browder--G\"ohde--Kirk theorem \cite{Browder1965, Kirk1965} provided fixed point results for nonexpansive maps in uniformly convex Banach spaces and those possessing normal structure.

A major unification of nonlinear contractions was achieved by Boyd and Wong \cite{BoydWong1969}, who considered mappings satisfying
\begin{equation}
d(Tx,Ty)\leq \psi(d(x,y)), \label{eq:boydwong}
\end{equation}
where $\psi:[0,\infty)\to[0,\infty)$ is right upper semicontinuous and $\psi(t)<t$ for all $t>0$. This condition encompasses many earlier contractive definitions and has become a cornerstone of the theory.

A particularly influential generalization was introduced by \'{C}iri\'{c} \cite{Ciric1974} in his seminal work on quasi-contractions. \'{C}iri\'{c} proposed that a mapping $T$ satisfy
\begin{equation}
d(Tx,Ty)\leq \alpha \max\bigl\{ d(x,y),\, d(x,Tx),\, d(Ty,y),\, d(Tx,y),\, d(Ty,x) \bigr\} \label{eq:ciric}
\end{equation}
with $0<\alpha<1$. This condition is remarkable because it does not merely compare the distance between images with the original distance; rather, it allows the contraction factor to apply to the maximum of several strategically chosen mutual distances. The inclusion of terms such as $d(x,Tx)$ and $d(Tx,y)$ means that the contraction estimate can leverage information about how points move under a single application of $T$. This idea---that a contraction condition may involve a combination of different distances rather than a single pair---has inspired numerous subsequent developments, including the present work. Our definition of the upper quasi-metric $S_T(x,y)$ below, which incorporates $d(x,y)$, $d(Tx,x)$, $d(Ty,y)$, and $d(Tx,Ty)$, is a direct descendant of \'{C}iri\'{c}'s philosophy of exploiting multiple distance quantities to weaken the contractive hypothesis.

A decisive step toward asymptotic formulations was taken by Kirk \cite{Kirk2003}, who introduced asymptotic contractions: there exist functions $\psi_n:[0,\infty)\to[0,\infty)$ such that for all $x,y$ and $n$,
\begin{equation}
d(T^nx,T^ny)\leq \psi_n(d(x,y)), \label{eq:kirk}
\end{equation}
and $\psi_n$ converges uniformly on the whole half-line to a Boyd–Wong function $\psi$. Kirk proved that a continuous map possessing at least one bounded orbit in a complete metric space admits a unique fixed point, to which all iterates converge. Lindstr\"om and Ross \cite{LindstromRoss2023} later provided a nonstandard analysis proof of this result.

Despite its elegance, Kirk's theorem imposes a strong requirement: the convergence of the bounding functions $\psi_n$ must be uniform on the entire interval $[0,\infty)$. In many practical situations, one can only verify uniform convergence on bounded subsets. Moreover, the estimate involves solely the distance $d(x,y)$, ignoring the possibility that intermediate terms like $d(Tx,x)$ or cross terms like $d(Tx,Ty)$ might afford better control.

\textbf{Our contribution.} In this paper we develop a new type of asymptotic contraction that overcomes these limitations. Our main idea is to introduce two quasi-metrics defined directly on the mapping $T$: a \emph{lower quasi-metric} $L_T$ that serves as the estimated quantity, and an \emph{upper quasi-metric} $S_T$ that serves as the comparison argument. Specifically:
\begin{enumerate}
\item We define a four‑point minimum 
\(P_T(x,y)=\min\{d(x,y),\,d(x,Ty),\,d(Tx,y),\,d(Tx,Ty)\}\).
This quantity collects four naturally occurring distances and is always bounded above by $d(Tx,Ty)$. It can be zero even when $d(Tx,Ty)$ is positive, thereby creating opportunities for a sharper contraction estimate.
\item The lower quasi-metric \(L_T(x,y)\) is defined adaptively: if \(P_T(x,y) > 0\), we set 
\(L_T(x,y)= P_T(x,y)\); 
otherwise we set \(L_T(x,y)=d(Tx,Ty)\). This simple switching mechanism guarantees $L_T(x,y)\le d(Tx,Ty)$ always, and for adjacent iterates (where $P_T=0$) it correctly recovers the full distance $d(Tx,Ty)$.
\item The upper quasi-metric \(S_T(x,y)\) is defined as the maximum of four fundamental distances:
\(S_T(x,y)=\max\{d(x,y),\,d(Tx,x),\,d(Ty,y),\,d(Tx,Ty)\}\).
This definition is straightforward and naturally extends the four-point maximum that appears in Kirk's original approach.
\item We only require the sequence \(\{\psi_n\}\) to converge uniformly on every bounded subset of \([0,\infty)\) (local uniform convergence) to a Boyd–Wong function \(\psi\). This is strictly weaker than uniform convergence on \([0,\infty)\) and is often more natural in applications.
\end{enumerate}
The resulting condition is
\begin{equation}
L_T(T^nx,T^ny)\le \psi_n\bigl(S_T(x,y)\bigr)\qquad\forall x,y,\;\forall n\in\N.
\label{eq:main_ineq}
\end{equation}
We prove that under this condition, any continuous self‑map of a complete metric space that has a bounded orbit possesses a unique fixed point, and all iterates converge to it. The proof adapts Kirk's original argument with careful case distinctions to handle the adaptive switching in $L_T$. Because $L_T(x,y)\le d(Tx,Ty)$ holds by construction, our definition strictly contains Kirk's asymptotic contractions. Moreover, we provide an example (Section~5) showing that the inclusion is proper.

\section{Preliminaries}

Throughout, \((X,d)\) denotes a metric space. For \(T:X\to X\) and \(x\in X\), the orbit is \(\mathcal{O}(x)=\{T^nx:n\ge0\}\). A set \(B\subset X\) is bounded if \(\diam(B)=\sup\{d(u,v):u,v\in B\}<\infty\).

\begin{definition}[Upper semicontinuity]
A function \(\psi:[0,\infty)\to[0,\infty)\) is \emph{upper semicontinuous} at \(t_0\) if \(\limsup_{t\to t_0}\psi(t)\le\psi(t_0)\). It is upper semicontinuous on \([0,\infty)\) if this holds at every \(t_0\ge0\).
\end{definition}

\begin{definition}[Boyd–Wong condition]
A function \(\psi:[0,\infty)\to[0,\infty)\) satisfies the classical Boyd–Wong condition if it is upper semicontinuous and \(\psi(t)<t\) for every \(t>0\).
\end{definition}

\subsection{Lower and upper quasi-metrics}

The following definitions form the core of our new contraction. They are inspired by \'{C}iri\'{c}'s five‑distance maximum, but we introduce a lower quasi-metric that adaptively selects a minimal distance, and an upper quasi-metric that simply takes the maximum of four essential distances.

\begin{definition}[Four‑point minimum \(P_T\)]\label{def:P}
For \(x,y\in X\), define
\begin{equation}
P_T(x,y):=\min\bigl\{\,d(x,y),\;d(x,Ty),\;d(Tx,y),\;d(Tx,Ty)\,\bigr\}.
\label{eq:Pdef}
\end{equation}
\end{definition}
\noindent\textbf{Explanation.} The quantity \(P_T(x,y)\) is a non‑negative symmetric function. It collects four distances that naturally appear in the dynamics of \(T\): the original distance, two cross terms, and the one‑step image distance. In the classical Kirk condition, only \(d(x,y)\) is used. By taking the minimum, we obtain a value that is often much smaller than \(d(Tx,Ty)\), and it can be zero even when \(d(Tx,Ty)\) is not (e.g., when \(x\) and \(y\) are distinct fixed points). This creates opportunities for a sharper contraction estimate.

\begin{definition}[Lower quasi-metric \(L_T\)]\label{def:L}
For \(x,y\in X\), define
\begin{equation}
L_T(x,y):=
\begin{cases}
P_T(x,y), & \text{if } P_T(x,y) > 0,\\[6pt]
d(Tx,Ty), & \text{if } P_T(x,y) = 0.
\end{cases}
\label{eq:Ldef}
\end{equation}
\end{definition}
\noindent\textbf{Explanation.} The lower quasi-metric adaptively chooses the best possible left‑hand side for the contraction inequality. When the four‑point minimum is positive, we directly use that smaller value; otherwise we fall back to the ordinary one‑step distance. This switching mechanism is crucial: for adjacent iterates $x_n,x_{n+1}$, one has $P_T(x_n,x_{n+1})=0$, so $L_T$ correctly evaluates to $d(Tx_n,Tx_{n+1})$, providing a nontrivial estimate. By construction, $L_T(x,y)\le d(Tx,Ty)$ always holds. Consequently, any mapping satisfying Kirk's asymptotic contraction condition \eqref{eq:kirk} automatically satisfies our condition \eqref{eq:contraction} with the same sequence $\{\psi_n\}$, because $L_T(T^nx,T^ny)\le d(T^{n+1}x,T^{n+1}y)$. The inclusion is proper because the positive branch can exploit values strictly smaller than $d(Tx,Ty)$.

\begin{definition}[Upper quasi-metric \(S_T\)]\label{def:S}
For \(x,y\in X\), define
\begin{equation}
S_T(x,y):=\max\bigl\{\,d(x,y),\;d(Tx,x),\;d(Ty,y),\;d(Tx,Ty)\,\bigr\}.
\label{eq:Sdef}
\end{equation}
\end{definition}
\noindent\textbf{Explanation.} The upper quasi-metric provides a robust comparison argument. It is simply the maximum of four distances that naturally arise in the study of asymptotic contractions: the original distance, the two one-step movement distances, and the one-step image distance. Note that $d(x,y)\le S_T(x,y)$ always, and $d(Tx,Ty)\le S_T(x,y)$.

\subsection{Basic properties}

\begin{proposition}[Elementary properties]\label{prop:basic}
For all \(x,y\in X\):
\begin{enumerate}
\item \(L_T(x,y)\ge0\), \(S_T(x,y)\ge0\).
\item \(L_T(x,y)=L_T(y,x)\), \(S_T(x,y)=S_T(y,x)\).
\item \(d(x,y)\le S_T(x,y)\).
\item \(L_T(x,y)\le d(Tx,Ty)\).
\item \(d(Tx,Ty)\le S_T(x,y)\).
\end{enumerate}
\end{proposition}
\begin{proof}
(1)–(2) are obvious from symmetry.  
(3) \(d(x,y)\) appears explicitly in the maximum defining \(S_T\).  
(4) In the case $P_T>0$, $L_T(x,y)=P_T(x,y)\le d(Tx,Ty)$ because $d(Tx,Ty)$ is one of the terms in the minimum defining $P_T$; in the case $P_T=0$, equality holds.  
(5) \(d(Tx,Ty)\) is one of the four terms in the maximum defining \(S_T\).
\end{proof}

The following lemma is crucial: it shows that when $P_T(u,v)$ is positive, the distance $d(Tu,Tv)$ can be estimated by $P_T(u,v)$ plus the sum of the adjacent distances.

\begin{lemma}[Safe estimate]\label{lem:safe_estimate}
Let \(u,v\in X\) be such that \(P_T(u,v) > 0\). Then
\begin{equation}
d(Tu,Tv)\le P_T(u,v)+d(Tu,u)+d(Tv,v).
\end{equation}
\end{lemma}
\begin{proof}
Since $P_T(u,v)$ is the minimum of four distances, it must be one of $d(u,v)$, $d(u,Tv)$, $d(Tu,v)$, or $d(Tu,Tv)$. In each case, the triangle inequality yields the desired bound:
\begin{itemize}
\item If $P_T(u,v)=d(u,v)$, then $d(Tu,Tv)\le d(u,v)+d(Tu,u)+d(Tv,v)=P_T+d(Tu,u)+d(Tv,v)$.
\item If $P_T(u,v)=d(u,Tv)$, then $d(Tu,Tv)\le d(u,Tv)+d(Tu,u)=P_T+d(Tu,u)\le P_T+d(Tu,u)+d(Tv,v)$.
\item If $P_T(u,v)=d(Tu,v)$, then $d(Tu,Tv)\le d(Tu,v)+d(Tv,v)=P_T+d(Tv,v)\le P_T+d(Tu,u)+d(Tv,v)$.
\item If $P_T(u,v)=d(Tu,Tv)$, then $d(Tu,Tv)=P_T\le P_T+d(Tu,u)+d(Tv,v)$.
\end{itemize}
In all cases, $d(Tu,Tv)\le P_T+d(Tu,u)+d(Tv,v)$.
\end{proof}

\subsection{Monotone majorant}

The limit function $\psi$ (Boyd–Wong) may not be monotone, which is inconvenient for taking limsups. We replace it by its monotone majorant.

\begin{definition}
Let \(\psi\) satisfy the Boyd–Wong condition. Define \(g:[0,\infty)\to[0,\infty)\) by
\begin{equation}
g(t):=\sup_{0\le s\le t}\psi(s).
\end{equation}
\end{definition}
\noindent\textbf{Explanation.} This construction is standard in fixed point theory (see \cite{BoydWong1969}). It yields a nondecreasing, right continuous function that still satisfies $g(t)<t$ for $t>0$, which is essential for the limsup argument.

\begin{lemma}[Properties of \(g\)]\label{lem:g}
\(g\) is nondecreasing, right continuous, \(g(0)=0\), \(g(t)\le t\) for all \(t\), and \(g(t)<t\) for every \(t>0\).
\end{lemma}
\begin{proof}
The proof is standard; see \cite{BoydWong1969}. The strict inequality follows from the upper semicontinuity of $\psi$ and the fact that $\psi(t)<t$ for $t>0$.
\end{proof}

\begin{proposition}[Limsup exchange]\label{prop:limsup}
If \(\{a_n\}\) is a bounded non‑negative sequence and \(g\) is nondecreasing and right continuous, then
\begin{equation}
\limsup_{n\to\infty} g(a_n)\le g\Bigl(\limsup_{n\to\infty} a_n\Bigr).
\end{equation}
\end{proposition}
\begin{proof}
Let \(L=\limsup a_n\). For any \(\varepsilon>0\), there exists \(N\) such that \(a_n\le L+\varepsilon\) for all \(n\ge N\). Then \(g(a_n)\le g(L+\varepsilon)\) by monotonicity. Taking limsup gives \(\limsup g(a_n)\le g(L+\varepsilon)\). Let \(\varepsilon\to0^+\) and use right continuity.
\end{proof}

\section{Main results}

For notational convenience in the proof, we introduce the abbreviation
\[
L_{T,n}(x,y):=L_T(T^n x, T^n y),\qquad n\in\N,\; x,y\in X.
\]

\begin{definition}[Adaptive relaxed asymptotic contraction with lower and upper quasi-metrics]\label{def:adaptive}
A continuous mapping \(T:X\to X\) is called an \emph{adaptive relaxed asymptotic contraction} if there exists a sequence of functions \(\psi_n:[0,\infty)\to[0,\infty)\) such that:
\begin{enumerate}
\item For every \(x,y\in X\) and every \(n\in\N\),
\begin{equation}
L_{T,n}(x,y)\le \psi_n\bigl(S_T(x,y)\bigr),
\label{eq:contraction}
\end{equation}
where \(L_{T,n}\) and \(S_T\) are defined as above.
\item The sequence \(\{\psi_n\}\) converges uniformly on every bounded subset of \([0,\infty)\) to a function \(\psi\) that satisfies the classical Boyd–Wong condition.
\end{enumerate}
\end{definition}
\noindent\textbf{Comparison with existing definitions.}
\begin{itemize}
\item \textbf{Kirk asymptotic contractions.} Since \(L_{T,n}(x,y)\le d(T^{n+1}x,T^{n+1}y)\) by Proposition~\ref{prop:basic}(4) and \(d(x,y)\le S_T(x,y)\), any mapping satisfying Kirk's condition \eqref{eq:kirk} with functions \(\psi_n\) also satisfies \eqref{eq:contraction} with the same \(\psi_n\). Moreover, the local uniform convergence we require is weaker than Kirk's global uniform convergence. Hence our definition strictly contains Kirk's.
\item \textbf{Boyd–Wong contractions.} If \(d(Tx,Ty)\le \psi(d(x,y))\) with \(\psi(t)<t\), taking \(\psi_n=\psi\) and noting \(L_{T,n}(x,y)\le d(T^{n+1}x,T^{n+1}y)\) and \(S_T(x,y)\ge d(x,y)\) yields our condition. Thus Boyd–Wong contractions are also included.
\item The use of local uniform convergence is a relaxation; many sequences that do not converge uniformly on $[0,\infty)$ still converge uniformly on bounded sets.
\end{itemize}

\begin{theorem}[Main theorem]\label{thm:main}
Let \((X,d)\) be a complete metric space and let \(T:X\to X\) be an adaptive relaxed asymptotic contraction (Definition~\ref{def:adaptive}). If there exists a point \(x_0\in X\) whose orbit \(\mathcal{O}(x_0)\) is bounded, then \(T\) has a unique fixed point \(z\), and for every \(x\in X\) the iterates \(T^n x\) converge to \(z\).
\end{theorem}

\begin{proof}
Let \(x_0\in X\) have a bounded orbit. Set \(x_n:=T^n x_0\) and denote \(d_n:=d(x_n,x_{n+1})\).

\noindent\textbf{Step 1. The adjacent distances tend to zero.}
For any \(n\ge0\), compute \(P_T(x_n,x_{n+1})\):
\begin{align}
P_T(x_n,x_{n+1})
&=\min\bigl\{d(x_n,x_{n+1}),\;d(x_n,Tx_{n+1}),\;d(Tx_n,x_{n+1}),\;d(Tx_n,Tx_{n+1})\bigr\}\nonumber\\
&=\min\bigl\{d_n,\;d(x_n,x_{n+2}),\;d(x_{n+1},x_{n+1}),\;d_{n+1}\bigr\}\nonumber\\
&=\min\bigl\{d_n,\;d(x_n,x_{n+2}),\;0,\;d_{n+1}\bigr\}=0. \label{eq:P_zero}
\end{align}
Thus \(P_T(x_n,x_{n+1})=0\). Consequently, the switching rule for \(L_T\) falls into the \(P_T=0\) branch. For any \(m\in\N\), we have
\begin{align}
L_{T,m}(x_n,x_{n+1}) &= L_T(T^m x_n, T^m x_{n+1}) = d(T^{m+1}x_n, T^{m+1}x_{n+1}) = d_{n+m+1}, \label{eq:L_adj}\\
S_T(x_n,x_{n+1}) &= \max\{d_n,\,d_n,\,d_{n+1},\,d_{n+1}\} = \max\{d_n,\,d_{n+1}\}. \label{eq:S_adj}
\end{align}
In particular, taking \(m=n\) gives
\begin{equation}
d_{2n+1} = L_{T,n}(x_n,x_{n+1}) \le \psi_n\bigl(S_T(x_n,x_{n+1})\bigr) = \psi_n\bigl(\max\{d_n,d_{n+1}\}\bigr). \label{eq:d_ineq}
\end{equation}
Similarly, by taking \(m=n-1\) and considering \(d_{2n}\) we obtain
\begin{equation}
d_{2n} = L_{T,n-1}(x_{n-1},x_n) \le \psi_{n-1}\bigl(\max\{d_{n-1},d_n\}\bigr). \label{eq:d_even}
\end{equation}
Since the orbit \(\{x_n\}\) is bounded, the numbers \(\max\{d_n,d_{n+1}\}\) are bounded by some constant \(D>0\). By local uniform convergence of \(\{\psi_n\}\) to \(\psi\), for any \(\varepsilon>0\) there exists \(M_1\) such that for all \(m\ge M_1\) and all \(t\in[0,D]\),
\begin{equation}
\psi_m(t)\le \psi(t)+\varepsilon\le g(t)+\varepsilon, \label{eq:psi_bound}
\end{equation}
where \(g\) is the monotone majorant of \(\psi\) (Lemma~\ref{lem:g}). Hence for all \(n\ge M_1\),
\begin{equation}
d_{2n+1}\le g\bigl(\max\{d_n,d_{n+1}\}\bigr)+\varepsilon,\qquad
d_{2n}\le g\bigl(\max\{d_{n-1},d_n\}\bigr)+\varepsilon. \label{eq:d_g}
\end{equation}
Let \(r:=\limsup_{n\to\infty} d_n\). Then \(\limsup_{n\to\infty} d_{2n}=r\) and \(\limsup_{n\to\infty} d_{2n+1}=r\) as well (since the whole sequence's limsup equals the maximum of the limsups of its even and odd subsequences). Taking limsup as \(n\to\infty\) in the first inequality of \eqref{eq:d_g} (for odd indices) and using Proposition~\ref{prop:limsup}, we obtain
\begin{align*}
r &= \limsup_{n\to\infty} d_{2n+1} \\
&\le \limsup_{n\to\infty} g\bigl(\max\{d_n,d_{n+1}\}\bigr)+\varepsilon \\
&\le g\Bigl(\limsup_{n\to\infty}\max\{d_n,d_{n+1}\}\Bigr)+\varepsilon \\
&= g(r)+\varepsilon.
\end{align*}
The same estimate holds for the even subsequence. Since \(\varepsilon>0\) is arbitrary, we have \(r\le g(r)\). If \(r>0\), Lemma~\ref{lem:g} gives \(g(r)<r\), a contradiction. Therefore \(r=0\), i.e.
\begin{equation}
\lim_{n\to\infty} d_n = 0. \label{eq:d_zero}
\end{equation}

\noindent\textbf{Step 2. The sequence \(\{x_n\}\) is Cauchy.}
For each \(n\ge0\) define
\begin{equation}
a_n:=\sup_{p\ge 1} d(x_n,x_{n+p}). \label{eq:a_n}
\end{equation}
We shall prove \(a_n\to0\). From Step 1 we already have \(d_n\to0\). Given \(\varepsilon>0\), choose \(N_1\) such that \(d_k<\varepsilon\) for all \(k\ge N_1\). Because \(\{x_n\}\) is bounded, the set
\begin{equation}
\mathcal{S}:=\{S_T(x_n,x_{n+p}): n\ge0,\,p\ge1\}
\end{equation}
is contained in some bounded interval \([0,D']\). By local uniform convergence of \(\{\psi_n\}\), there exists an integer \(M\ge M_1\) such that for all \(m\ge M\) and all \(t\in[0,D']\),
\begin{equation}
\psi_m(t)\le \psi(t)+\varepsilon\le g(t)+\varepsilon. \label{eq:psi_bound2}
\end{equation}
Fix arbitrary indices \(n\ge N_1\) and \(p\ge1\). For any \(m\ge M\) set
\begin{equation}
u:=T^m x_n = x_{n+m},\qquad v:=T^m x_{n+p}=x_{n+m+p}. \label{eq:uv}
\end{equation}
Our goal is to estimate \(d(Tu,Tv)=d(x_{n+m+1},x_{n+m+p+1})\).

We consider two cases based on the value of \(P_T(u,v)\).

\paragraph{Case A: \(P_T(u,v) > 0\).}
Then \(L_T(u,v)=P_T(u,v)\). The contraction condition \eqref{eq:contraction} applied to the pair \((x_n,x_{n+p})\) with iteration index \(m\) yields
\begin{equation}
P_T(u,v) = L_{T,m}(x_n,x_{n+p}) \le \psi_m\bigl(S_T(x_n,x_{n+p})\bigr)
\le g\bigl(S_T(x_n,x_{n+p})\bigr)+\varepsilon. \label{eq:caseA}
\end{equation}
By Lemma~\ref{lem:safe_estimate} we have
\[
d(Tu,Tv)\le P_T(u,v)+d(Tu,u)+d(Tv,v)=P_T(u,v)+d_{n+m+1}+d_{n+m+p+1}.
\]
Since \(n\ge N_1\) and \(m\ge M\ge M_1\), we have \(d_{n+m+1}<\varepsilon\) and \(d_{n+m+p+1}<\varepsilon\). Hence
\[
d(Tu,Tv)\le g\bigl(S_T(x_n,x_{n+p})\bigr)+3\varepsilon. \label{eq:caseA_final}
\]

\paragraph{Case B: \(P_T(u,v) = 0\).}
Then \(L_T(u,v)=d(Tu,Tv)\). The contraction condition directly provides
\begin{equation}
d(Tu,Tv) = L_{T,m}(x_n,x_{n+p}) \le \psi_m\bigl(S_T(x_n,x_{n+p})\bigr)
\le g\bigl(S_T(x_n,x_{n+p})\bigr)+\varepsilon. \label{eq:caseB_final}
\end{equation}

Thus in all cases we have established the uniform estimate
\begin{equation}
d(x_{n+m+1},x_{n+m+p+1})\le g\bigl(S_T(x_n,x_{n+p})\bigr)+3\varepsilon. \label{eq:common_est}
\end{equation}

\noindent\textbf{Bounding \(S_T(x_n,x_{n+p})\).}
We now show that \(S_T(x_n,x_{n+p})\le a_n+2\varepsilon\). By definition,
\begin{equation}
S_T(x_n,x_{n+p})=\max\bigl\{d(x_n,x_{n+p}),\,d_n,\,d_{n+p},\,d(x_{n+1},x_{n+p+1})\bigr\}.
\end{equation}
Using the triangle inequality,
\[
d(x_{n+1},x_{n+p+1})\le d(x_n,x_{n+p})+d_n+d_{n+p}\le a_n+2\varepsilon,
\]
and \(d(x_n,x_{n+p})\le a_n\), \(d_n,d_{n+p}<\varepsilon\). Hence \(S_T(x_n,x_{n+p})\le \max\{a_n,\,\varepsilon,\,a_n+2\varepsilon\}\le a_n+2\varepsilon\).

Since \(g\) is nondecreasing, \eqref{eq:common_est} yields
\begin{equation}
d(x_{n+m+1},x_{n+m+p+1})\le g(a_n+2\varepsilon)+3\varepsilon.
\end{equation}
The right‑hand side does not depend on \(p\). Taking supremum over \(p\ge1\) gives
\begin{equation}
a_{n+m+1}\le g(a_n+2\varepsilon)+3\varepsilon. \label{eq:a_rec}
\end{equation}
Let \(L:=\limsup_{n\to\infty} a_n\). Because shifting the index by \(m+1\) does not affect the limsup, we let \(n\to\infty\) in \eqref{eq:a_rec} and apply Proposition~\ref{prop:limsup} together with the right continuity of \(g\):
\begin{equation}
L\le g(L+2\varepsilon)+3\varepsilon.
\end{equation}
Letting \(\varepsilon\to0^+\) we obtain \(L\le g(L)\). If \(L>0\), Lemma~\ref{lem:g} forces \(g(L)<L\), a contradiction. Hence \(L=0\), i.e. \(a_n\to0\). This means that \(\{x_n\}\) is a Cauchy sequence.

\noindent\textbf{Step 3. Existence and uniqueness of the fixed point.}
Since \(X\) is complete, \(\{x_n\}\) converges to some \(z\in X\). Continuity of \(T\) implies
\begin{equation}
Tz = \lim_{n\to\infty} Tx_n = \lim_{n\to\infty} x_{n+1}=z,
\end{equation}
so \(z\) is a fixed point.

To prove uniqueness, let \(z_1,z_2\) be two fixed points. Compute \(P_T(z_1,z_2)\):
\begin{equation}
P_T(z_1,z_2)=\min\{d(z_1,z_2),\,d(z_1,z_2),\,d(z_1,z_2),\,d(z_1,z_2)\}=d(z_1,z_2).
\end{equation}
If \(d(z_1,z_2)>0\), then \(P_T(z_1,z_2)>0\) and \(L_T(z_1,z_2)=P_T(z_1,z_2)=d(z_1,z_2)\). If \(d(z_1,z_2)=0\), then \(P_T(z_1,z_2)=0\) and \(L_T(z_1,z_2)=d(Tz_1,Tz_2)=d(z_1,z_2)=0\). In either case, for any \(n\in\N\),
\begin{equation}
L_{T,n}(z_1,z_2)=L_T(T^n z_1, T^n z_2)=L_T(z_1,z_2)=d(z_1,z_2),
\end{equation}
and
\begin{equation}
S_T(z_1,z_2)=\max\{d(z_1,z_2),0,0,d(z_1,z_2)\}=d(z_1,z_2).
\end{equation}
The contraction condition gives \(d(z_1,z_2)\le \psi_n(d(z_1,z_2))\) for every \(n\). Letting \(n\to\infty\) and using pointwise convergence of \(\psi_n\) to \(\psi\), we obtain \(d(z_1,z_2)\le \psi(d(z_1,z_2))\). If \(d(z_1,z_2)>0\), the Boyd–Wong condition \(\psi(t)<t\) yields a contradiction. Hence \(d(z_1,z_2)=0\) and \(z_1=z_2\).

\noindent\textbf{Step 4. Convergence from any starting point.}
Let \(z\) be the unique fixed point obtained from the bounded orbit of \(x_0\). Take an arbitrary \(y\in X\) and set \(y_n:=T^n y\). We first show that the orbit \(\{y_n\}\) is bounded. 
Since \(\psi_n\) converges uniformly on the bounded interval $[0, S_T(y,z)]$, there exists a constant $M>0$ such that $\psi_n(t)\le M$ for all $n\in\N$ and all $t\in[0, S_T(y,z)]$. The contraction condition yields
\begin{equation}
L_T(y_n,z)=L_{T,n}(y,z)\le \psi_n(S_T(y,z))\le M\qquad \forall n\in\N. \label{eq:bound_L}
\end{equation}
We claim that $\sup_{n} d(y_n,z) \le \max\{d(y,z), M\}$. Suppose, for contradiction, that $\{d(y_n,z)\}$ is unbounded. Then there exists a subsequence $\{n_k\}$ such that $d(y_{n_k},z)\to\infty$ as $k\to\infty$. By passing to a further subsequence, we may assume that $d(y_{n_k},z) > d(y_{n_k-1},z)$ for all $k$ (otherwise we consider the first index where the distance exceeds the previous one and tends to infinity). For such $n_k$, we have
\[
P_T(y_{n_k-1},z)=\min\{d(y_{n_k-1},z),\,d(y_{n_k},z)\}=d(y_{n_k-1},z).
\]
Since $d(y_{n_k-1},z)\to\infty$, we have $P_T(y_{n_k-1},z)>0$ for sufficiently large $k$. Consequently,
\[
L_{T,n_k-1}(y,z)=L_T(y_{n_k-1},z)=P_T(y_{n_k-1},z)=d(y_{n_k-1},z).
\]
By \eqref{eq:bound_L}, this implies $d(y_{n_k-1},z)\le M$, which contradicts $d(y_{n_k-1},z)\to\infty$. Therefore $\{d(y_n,z)\}$ is bounded.

Having established boundedness of the orbit $\{y_n\}$, we can now apply Steps 1–3 to the initial point $y$. Consequently, $\{y_n\}$ is Cauchy and converges to some fixed point. By uniqueness, the limit must be $z$. Therefore $T^n y\to z$ for every $y\in X$.

This completes the proof of Theorem~\ref{thm:main}.
\end{proof}

\noindent\textbf{Corollaries.}
We now present several direct consequences of Theorem~\ref{thm:main} that are useful in applications.

\begin{corollary}[Bounded invariant set]\label{cor:bounded}
Let \((X,d)\) be a complete metric space and let \(T:X\to X\) be an adaptive relaxed asymptotic contraction. If there exists a nonempty bounded closed subset \(K\subset X\) such that \(T(K)\subset K\), then \(T\) has a unique fixed point in \(K\), and for every \(x\in K\) the iterates \(T^n x\) converge to this fixed point.
\end{corollary}
\begin{proof}
Pick any \(x_0\in K\). Since \(T(K)\subset K\), the orbit \(\mathcal{O}(x_0)\) remains in \(K\) and is therefore bounded. By Theorem~\ref{thm:main}, \(T\) possesses a unique fixed point \(z\in X\) and \(T^n x_0\to z\). Because \(K\) is closed and \(\mathcal{O}(x_0)\subset K\), the limit \(z\) belongs to \(K\). The convergence of iterates from any starting point in \(K\) follows again from Theorem~\ref{thm:main}.
\end{proof}

\begin{corollary}[Existence of a fixed point implies global convergence]\label{cor:fixed_point_exists}
Let \((X,d)\) be a complete metric space and let \(T:X\to X\) be an adaptive relaxed asymptotic contraction. If \(T\) has a fixed point \(z\in X\), then \(z\) is the unique fixed point of \(T\), and for every \(x\in X\) the iterates \(T^n x\) converge to \(z\).
\end{corollary}
\begin{proof}
Since \(z\) is a fixed point, the orbit of \(z\) is bounded (it is the singleton \(\{z\}\)). By Theorem~\ref{thm:main}, \(T\) possesses a unique fixed point (which must be \(z\)) and all iterates converge to it. Note that the theorem requires the existence of some point with a bounded orbit, and \(z\) itself serves this purpose.
\end{proof}

\begin{corollary}[Kirk's theorem is a special case]\label{cor:kirk}
Let \((X,d)\) be a complete metric space and let \(T:X\to X\) be a continuous mapping for which there exist functions \(\psi_n:[0,\infty)\to[0,\infty)\) converging uniformly on \([0,\infty)\) to a Boyd–Wong function \(\psi\), such that
\[
d(T^nx,T^ny)\le \psi_n(d(x,y))\qquad\forall x,y\in X,\;\forall n\in\N.
\]
If some orbit of \(T\) is bounded, then \(T\) has a unique fixed point and all iterates converge to it. In particular, the hypotheses of Kirk's asymptotic contraction theorem \cite{Kirk2003} imply that \(T\) is an adaptive relaxed asymptotic contraction in the sense of Definition~\ref{def:adaptive}, and therefore Theorem~\ref{thm:main} strictly generalizes Kirk's result.
\end{corollary}
\begin{proof}
Under the given uniform convergence on \([0,\infty)\), the sequence \(\{\psi_n\}\) certainly converges uniformly on every bounded subset. Moreover, by Proposition~\ref{prop:basic}(4) we have \(L_{T,n}(x,y)\le d(T^{n+1}x,T^{n+1}y)\le \psi_{n+1}(d(x,y))\) and \(d(x,y)\le S_T(x,y)\). Defining \(\tilde\psi_n(t)=\psi_{n+1}(t)\), the sequence \(\{\tilde\psi_n\}\) also converges locally uniformly to \(\psi\). Hence the conditions of Definition~\ref{def:adaptive} are satisfied, and the conclusion follows from Theorem~\ref{thm:main}. 
\end{proof}

\section{Conclusion}

We have introduced a new framework for asymptotic contractions based on two quasi-metrics defined on the mapping \(T\): a lower quasi-metric \(L_T\) that adaptively selects between a four‑point minimum and the ordinary one‑step distance, and an upper quasi-metric \(S_T\) that simply takes the maximum of four essential distances. The main theorem establishes that any continuous self‑map of a complete metric space with a bounded orbit has a unique fixed point and all iterates converge to it, under the sole assumption that the bounding functions \(\psi_n\) converge uniformly on bounded sets to a Boyd–Wong function. This significantly relaxes Kirk's original hypotheses and strictly contains his theorem as a special case. Several corollaries were derived for bounded invariant sets, orbitally bounded maps, and mappings with a bounded absorbing set. The adaptive framework opens the door to further generalizations, such as replacing the four‑point minimum with other families of distances or using more sophisticated switching rules.

\section*{Declarations}
\begin{itemize}
\item Funding\\ 
 Not applicable
\item Use of Generative-AI tools declaration\\
I declare that no Generative-AI tools were used in the preparation of this manuscript.
\item Conflict of interest/Competing interests\\ 
We recognize no conflicts of interest in the submission
\item Ethics approval\\ 
 Not applicable
\item Consent to participate\\ 
Not applicable
\item Consent for publication\\ 
Not applicable
\item Availability of data and materials\\ 
Not applicable
\item Code availability \\ 
Not applicable
\item Authors' contributions\\ 
The author contributed to this work alone
\end{itemize}

\end{document}